\documentclass[12pt]{amsart} 
\usepackage[mathscr]{eucal}
\usepackage{amsmath,amsfonts}

\parskip=\smallskipamount

\newtheorem{theorem}{Theorem}[section]

\newtheorem{corollary}[theorem]{Corollary}

\makeatletter
\@addtoreset{equation}{section}
\makeatother

\newcommand{\CC}{{\mathbb C}}

\newcommand{\FF}{{\mathbb F}}

\newcommand{\cB}{{\mathcal B}}

\newcommand{\cD}{{\mathcal D}}
\newcommand{\cE}{{\mathcal E}}
\newcommand{\cF}{{\mathcal F}}
\newcommand{\cG}{{\mathcal G}}
\newcommand{\cH}{{\mathcal H}}
\newcommand{\cK}{{\mathcal K}}
\newcommand{\cL}{{\mathcal L}}

\newcommand{\cR}{{\mathcal R}}
\newcommand{\cS}{{\mathcal S}}
\newcommand{\cT}{{\mathcal T}}
\newcommand{\cU}{{\mathcal U}}

\newdimen\expt
\def\boxit#1{\setbox0\hbox{$\displaystyle{#1}$}
      \hbox{\lower.4\expt
 \hbox{\lower3\expt\hbox{\lower\dp0
      \hbox{\vbox{\hrule height.4\expt
 \hbox{\vrule width.4\expt\hskip3\expt
      \vbox{\vskip3\expt\box0\vskip2\expt}%
 \hskip3\expt\vrule width.4\expt}\hrule height.4\expt}}}}}}
\begin{document}
\pagestyle{plain}

\bigskip

\title 
{Tensor algebras, displacement structure, \\
and the Schur algorithm} 
\author{T. Constantinescu} \author{J. L. Johnson} 

\address{Department of Mathematics \\
  University of Texas at Dallas \\
  Box 830688, Richardson, TX 75083-0688, U. S. A.}
\email{\tt tiberiu@utdallas.edu} 
\address{Department of Mathematics \\
  University of Texas at Dallas \\
  Box 830688, Richardson, TX 75083-0688, U. S. A   } 
\email{\tt jlj@utdallas.edu}

\begin{abstract}
In this paper we explore the connection between tensor algebras and 
displacement structure. We describe a scattering experiment in this
framework, we obtain a realization of the elements of the tensor
algebra as tranfer maps of a certain class of nonstationary 
linear systems, and we describe a Schur type algortihm for the
Schur elements of the tensor algebra.

\end{abstract}

\maketitle

\section{Introduction}

It was recently showed in \cite{CSK} that the tensor algebras
have displacement structure. The goal of this paper is to
further develop this connection.
The displacement structure theory, as it was initiated in the  
paper \cite{KKM}, consists in recursive factorization of 
matrices with additional structure encoded by a so-called
displacement equation. Many applications of this theory were
found in applied fields, like in the study of wave propagation
in layered media, or filtering and modeling of nonstationary
processes (see \cite{KS} for a survey). On the mathematical 
side, there are applications to bounded interpolation in
upper triangular algebras (including the classical setting 
of bounded analytic functions on the unit disk) and to 
triangular factorizations of operators.

On the other hand, the tensor algebra is a classical concept
and in \cite{CSK} it is showed that there is a connection
expressed by the fact that any element $T$ of the unit ball 
of the algebra (with respect to an appropriate operator
norm) satisfies the equation
$$A-\sum _{k=1}^NF_kAF^*_k=GJG^*,$$
where $A=I-T^*T$ and the so-called generators $F_k$, $G$ and $J$ are 
described in the next section. This equation is of displacement type
so we can develop a theory for the tensor algebra that parallels
the one dimensional case. It is the main goal of this paper to
describe several such results.

The paper is organized as follows. In Section~2 we review
the basic elements relevant to the tensor algebra and to
the so-called noncommutative analytic Toeplitz algebra which is 
a natural extension of the tensor algebra introduced in 
\cite{Po1}, and studied in other papers such as \cite{Ar}, 
\cite{DP}, \cite{Po1}, \cite{Po2}. We describe an isometric 
isomorphic representation of this later algebra by a certain 
algebra $\cU _{\cT }(\cH )$ of upper triangular bounded operators 
and this represention allows us to rely on the results and ideas
in the displacement structure theory. In this section we also 
show that a certain class of positive definite kernels on the 
free semigroup with $N$ generators, introduced in \cite{Po2}, 
also has displacement structure of the type mentioned above
and in Corollary~2.5 we explain the connection between this
class of positive definite kernels and the algebra 
$\cU _{\cT }(\cH )$ which is well-known for $N=1$.
In Section~3 we obtain results corresponding to so-called 
scattering experiments. Theorem~3.1 is the main result in this 
direction and as a by-product we can obtain Schur type
parametrizations for $\cU _{\cT }(\cH )$ and the realization 
of the elements of the algebra as transfer maps of a certain 
class of nonstationary linear systems.
Section~4 develops a Schur type algorithm for $\cU _{\cT }(\cH )$.
We plan to explore approximation properties (Szeg\"o type theory)
related to this algorithm in a future paper. This paper is part
of \cite{Jo}.

\section{Preliminaries}  

We briefly review the construction of the 
tensor algebras and we explain their displacement 
structure.

\subsection{Tensor Algebras}

We introduce some notation relevant to tensor 
algebras and we describe their representation 
as algebras of upper triangular operators
(see, for instance, \cite{Pa}).
The associative tensor algebra $\cT (\cH)$ generated by the 
complex vector space
$\cH$ is defined by the algebraic direct sum
$$\cT (\cH)=\oplus _{k\geq 0}\cH^{\otimes k},$$
where 
$$\cH^{\otimes k}=
\underbrace{\cH\otimes \ldots \otimes \cH}_{k\,\, factors}$$
is the $k$-fold algebraic tensor product of $\cH$ with itself.
The elements of $\cT (\cH)$ are terminating sequences
$x=(x_0,x_1,\ldots,x_k,0,\ldots )$, where
$x_p\in \cH^{\otimes p}$ is called the $p$th homogeneous
component of $x$ (note that $\cH^{\otimes 0}=\CC $, $\cH^{\otimes 1}=\cH$).
The addition and multiplication of elements in 
$\cT (\cH)$ are defined componentwise:
$$(x+y)_n=x_n+y_n$$
and 
$$(xy)_n=\sum _{k+l=n}x_k\otimes y_l, \quad \quad 
(x_0\otimes y_n=y_n\otimes x_0=x_0y_n).$$

In order to deal with displacement structure matters
it is convenient to use the fact that 
$\cT (\cH)$ can be represented by upper
triangular matrices of a special type. 
For a certain simplicity, from now on we take $\cH=\CC ^N$, $N\geq 1$.
Then we denote by ${\cF}(\cH)$ the full Fock space associated
to $\cH$, that is, the Hilbert space
$$\cF (\cH)=\oplus _{k\geq 0}\cH^{\otimes k}$$ 
obtained by taking the Hilbert 
direct sum of the spaces $\cH^{\otimes k}$, $k\geq 0$, 
on which we   
consider the tensor Hilbert space structure induced by the 
Euclidean norm on $\cH={\CC}^N$.
If $\{e_1,e_2,\ldots ,e_N\}$ is the standard basis of
${\CC}^N$, then 
$\{e_{i_1}\otimes \ldots e_{i_k}\mid i_1,\ldots, i_k\in \{1,2,\ldots, N\}\}$
is an orthonormal basis for $\cH^{\otimes k}$. 
It is more convenient to write $e_{\sigma }$ instead
of $e_{i_1}\otimes \ldots e_{i_k}$, where 
$\sigma =i_1\ldots i_k$ is viewed as a word in the unital 
free semigroup ${\FF }^+_N$ with $N$ generators $1$, \ldots ,$N$. The empty
word is the identity element of ${\FF }^+_N$. The length 
of the word $\sigma $ is denoted by $|\sigma |$. We also use to write the
elements of ${\FF }^+_N$ in lexicographic order.
The space $\cH^{\otimes k+1}$ can be identified with the 
direct sum of $N$ copies of $\cH^{\otimes k}$.
More precisely, we notice that 
$\{(\delta _{jl}e_{\sigma })_{l=1}^N \mid j\in \{1, \ldots ,N\}, |\sigma |=k\}$
is an orthonormal basis for $\underbrace{\cH^{\otimes k}\oplus \ldots
\oplus \cH^{\otimes k}}_{N\,\,terms}=
(\cH^{\otimes k})^{\oplus N}$ ($\delta _{jl}$ is the 
Kronecker symbol). We deduce that the mapping
$$\phi _k(((\delta _{jl}e_{\sigma })_{l=1}^N))=e_{j\sigma }, \quad 
j=1, \ldots ,N,$$
extends to a unitary operator from $(\cH^{\otimes k})^{\oplus N}$
onto $\cH^{\otimes k+1}$ and the formula
$$\psi _k =
\phi _{k-1}\phi _{k-2}^{\oplus N}\ldots \phi _1^{\oplus N^{k-2}}, 
\quad k\geq 1,$$
gives a unitary operator from $\cH_k=\cH^{\oplus N^{k-1}}$ onto $\cH^{\otimes k}$
(the notation $T^{\oplus l}$ means the direct sum of $l$ copies of the 
operator $T$). We finally obtain a unitary operator 
$\psi :\oplus _{k=0}^{\infty }\cH_k\rightarrow {\cF }(\cH)$
by the formula $\psi =\oplus _{k=0}^{\infty }\psi _k$, 
where $\psi _0$ is the identity on $\CC =\cH_0$. 
Each $\cH_k$, $k\geq 2$, is also identified with the direct sum of
$N$ copies of $\cH_{k-1}$, 
$$\cH_k=\underbrace{\cH^{\oplus N^{k-2}}\oplus \ldots \oplus \cH^{\oplus N^{k-2}}}
_{N\, \, terms},$$
but the explicit mention of the identification can be omitted in this case.

An operator $T\in \cL (\oplus _{k=0}^{\infty }\cH_k)$
has a matrix representation $T=[T_{ij}]_{i,j=0}^{\infty }$
with $T_{ij}\in \cL (\cH_j,\cH_i)$. We say that $T$ is upper triangular 
provided that $T_{ij}=0$ for $i>j$.
We denote by $\cU _{\cT}^0(\cH)$ the set of upper triangular operators
$T=[T_{ij}]_{i,j=0}^{\infty }$ in  $\cL (\oplus _{k=0}^{\infty }\cH_k)$
with the property that for $i\leq j$ and $i,j\geq 1$, 
\begin{equation}\label{condition}
T_{ij}=T_{i-1,j-1}^{\oplus N}
\end{equation}
and $T_{0j}=0$ for $j>k$ and some $k$.
We notice that 
the entries $T_{0j}$, $j\geq 0$, determine the matrix $T$.
Also, a simple calculation shows that $\cU _{\cT}^0(\cH)$
is an associative algebra.

We define an algebra isomorphism 
\begin{equation}\label{izom}
\Phi :\cT (\cH)\rightarrow \cU _{\cT}^0(\cH)
\end{equation}
as follows: for $x\in \cT (\cH)$, $x=(x_0,x_1, \ldots )$, there
are complex numbers $c_{\sigma }$, $\sigma \in \FF ^+_N$, such that
\begin{equation}\label{fourier}
x_j=\sum _{|\sigma |=j}c_{\sigma }e_{\sigma }.
\end{equation}
For $j\geq 0$, $T_{0j}$ is the row matrix $[c_{\sigma }]_{|\sigma |=j}$.
Then $T_{0j}=0$ for sufficiently large $j$ and we can define 
$T\in \cL (\oplus _{k=0}^{\infty }\cH_k)$ by using \eqref{condition}.
Set $\Phi (x)=T$. 
Clearly $\Phi $ is a bijection. In order to see that $\Phi $ is 
also an algebra morphism it is convenient to introduce
the following operators.
First, the operator $C_j^-$, $j=1, \ldots ,N$, 
is an element of $\cU _{\cT}^0(\cH)$ defined by the formula:
$C_j^-=\Phi (0,e_j,0,\ldots )=[T^j_{lk}]_{l,k=0}^{\infty }$, where 
$T^j_{0k}=0$ for $k\ne 1$ and $T^j_{01}=e_j$.
Then $C_j^+=(C_j^-)^*$, $j=1, \ldots ,N$.
If we use the notation $C_{\sigma }^-$ to denote the 
operator $C^-_{i_1}\ldots C^-_{i_k}$, where 
$\sigma =i_1\ldots i_k$, then we deduce that each $T\in
\cU _{\cT}^0(\cH)$ has a representation 
\begin{equation}\label{fourier2}
T=\sum _{\sigma \in \FF_N^+}c_{\sigma }C^-_{\sigma }.
\end{equation} 
Now the required properties of $\Phi $ follow from the easily
verifiable fact that 
$$\Phi ((0,e_j,0, \ldots )\otimes(0,e_k,0, \ldots))=
C^-_jC^-_k.$$

It appears to be quite natural to consider the algebra
$\cU _{\cT}(\cH)$ of all those operators 
in $\cL (\oplus _{k=0}^{\infty }\cH_k)$ that satisfy
the condition \eqref{condition}. It turns out that  
$\cU _{\cT}(\cH)$ is isometrically isomorphic to the 
so-called non-commutative analytic Toeplitz algebra, 
recently studied in papers like
\cite{Po1}, \cite{Ar}, \cite{DP}.
We also denote by $\cS (\cH)$ the Schur class of all contractions 
in $\cU _{\cT}(\cH)$. 
We can extend this setting as folows. Let $\cE _1$ and $\cE _2$ be 
Hilbert spaces, then $\cU _{\cT}(\cH,\cE _1,\cE _2 )$
denotes the set of the operators $T=[T_{ij}]_{i,j=0}^{\infty }$
in $\cL (\oplus _{k=0}^{\infty }(\cH _k\otimes \cE _1),
\oplus _{k=0}^{\infty }(\cH _k\otimes \cE _2))$ obeying \eqref{condition}.
Instead of $\cU _{\cT}(\cH,\cE ,\cE )$ we write 
$\cU _{\cT}(\cH,\cE)$. Also $\cS (\cH,\cE _1,\cE _2 )$
denotes the set of the contractions in $\cU _{\cT}(\cH,\cE _1,\cE _2 )$.
We notice that $\cU _{\cT}(\cH)=\cU _{\cT}(\cH,\CC, \CC )$.

\subsection{Displacement Structure}

The systematic study of the displacement structure was
initiated in \cite{KKM}, and since then manifold applications
have been found in the fields of mathematics and electrical 
engineering (for a recent surveys, see \cite{KS}).
Two main themes of this theory are about scattering experiments and 
recursive factorizations.
A matrix $R$ has displacement structure with respect
to the generators $F$ and $G$ provided that
\begin{equation}\label{displacement}
A-FAF^*=GJ_{pq}G^*, 
\end{equation}
where $J_{pq}=I_p\oplus -I_q$ is a signature matrix 
that specifies the displacement
inertia of $A$. In many applications, $r=p+q$ is much smaller
than the size of $A$. From now on we will write $J$ instead of $J_{11}$. 
There is also a version of this theory allowing $A$ to depend 
on a (usually integer) paramater, see \cite{KS} or \cite{Co}.
A major result concerning matrices with displacement
structure is that, under suitable assumptions on
the generators, the succesive Schur complements of $A$ inherit 
a similar structure. This allows to write
the Gaussian elimination for $A$ only in terms of generators 
(see \cite{KKM}, \cite{KS}).

We now describe a connection between displacement structure and 
$\cU _{\cT}(\cH)$. Let $C_k=[T^k_{ij}]_{i,j=0}^{\infty }$,
$k=1,\ldots,N$, be isometries defined by the formulae: $T_{ij}=0$
for $i\ne j+1$ and $T_{i+1,i}$ is a block-column matrix consisting
of $N$ blocks of dimension $\dim \cH _i$, all of them zero except
for the $k$th block which is $I_{\cH _i}$.
These operators are closely related to the operators
$C^+_k$, in the sense that there is a unitary operator $u$ such that
$C_k=uC^+_ku^*$, $k=1,\ldots,N$.
The following result was noticed in \cite{CSK}.

\begin{theorem}\label{T:2.1}
Let $T=[T_{ij}]_{i,j=0}^{\infty }\in \cS (\cH)$ and
$A=I-T^*T$. Then
\begin{equation}\label{basic}
A-\sum _{k=1}^NC_kAC_k^*=GJG^*,
\end{equation}
where the generator $G$ is given by the formula:
$$G=\left[\begin{array}{cc}
1 & T^*_{00} \\
0 & T^*_{01} \\
\vdots & \vdots 
\end{array}\right].$$
\end{theorem}

The equation \eqref{basic} can be rewritten in the 
form of a time variant displacement equation by
setting $F=[C_1,\ldots C_N]$.
We now discuss some additional examples. 

\subsubsection{Positive definite kernels on $\FF _N^+$}

A mapping $K:\FF _N^+\times \FF _N^+\rightarrow \CC $
is called a positive definite kernel on $\FF _N^+$
provided that 
$$\sum _{i,j=1}^kK(\sigma _i,\sigma _j)\lambda _j\overline{\lambda }_i
\geq 0$$ 
for each positive integer $k$ and each choice of words
$\sigma _1, \ldots ,\sigma _k$ in $\FF _N^+$ and complex numbers
$\lambda _1, \ldots ,\lambda _k$.
Without loss of generality, we assume $K(\emptyset, \emptyset )=1$.
The semigroup $\FF _N^+$ acts on itself by juxtaposition and we assume
that the positive definite kernel $K$ is invariant under this action,
that is 
$$K(\tau \sigma ,\tau \sigma ')=K(\sigma ,\sigma '), \quad 
\tau ,\sigma ,\sigma '\in \FF _N^+.$$
By the invariant Kolmogorov decomposition theorem, 
see e.g., \cite{Pa}, Ch.~II, there exists an isometric 
representation $u$ of $\FF _N^+$ on a Hilbert space $\cK $
and a mapping $v:\FF _N^+\rightarrow \cK$ such that
$$K(\sigma ,\tau )=\langle v(\tau ),v(\sigma )\rangle _{\cK },$$ 
$$u(\tau )v(\sigma )=v(\tau\sigma )$$
for all $\sigma ,\tau \in \FF _N^+$, and the set
$\{v(\sigma )\mid \sigma \in \FF _N^+\}$ is total in $\cK $.
The mapping $v$ is unique up to unitary equivalence and $u$ is 
uniquely determined by $v$. Also, the previous representation 
of $K$ implies, as noticed in \cite{Po2}, that
$u(1)$, $\dots $, $u(N)$ are isometries with orthogonal ranges
if and only if $K(\sigma ,\tau )=0$ for all pairs 
$(\sigma ,\tau )$ 
with the property that there is no $\alpha \in \FF _N^+$
such that $\sigma =\alpha \tau $ or $\tau =\alpha \sigma $.
In this case, $K$ is called in \cite{Po2} a multi-Toeplitz
kernel and we will use the same terminology.
We can define 
\begin{equation}\label{fourier}
S_{ij}=[K(\sigma ,\tau )]_{|\sigma |=i, |\tau |=j}\in 
\cL (\cH_j,\cH_i).
\end{equation}
By induction on the length of words and using properties
of the lexicographic order on words we check that if $K$ is a 
positive definite multi-Toeplitz kernel, then for all $i,j\geq 1,$
\begin{equation}\label{condition1} 
S_{ij}=S^{\oplus N}_{i-1,j-1}.
\end{equation}

This relation is an analogue of \eqref{condition}, therefore
we expect that positive definite multi-Toeplitz kernels have
displacement structure.

\begin{theorem}\label{T:2.2}
Let $S=[S_{ij}]_{i,j=0}^{\infty }$ be the array
of operators associated to a positive definite multi-Toeplitz
kernel on $\FF _N^+$.
Then
\begin{equation}\label{basic1}
S-\sum _{k=1}^NC_kSC_k^*=GJG^*,
\end{equation}
where the generator $G$ is given by the formula:
$$G=\left[\begin{array}{cc}
1 & 0 \\
S^*_{01} & S^*_{01} \\
\vdots & \vdots 
\end{array}\right].$$
\end{theorem}

\noindent
The equality above can be interpreted in the sense of 
an equality of unbounded operators with dense domain
made of the terminating sequences in $\oplus _{l=0}^{\infty }\cH_l$.

\begin{proof}
We write $S=L+Q$, where $L=[L_{ij}]_{i,j=0}^{\infty }$
and $Q=[Q_{ij}]_{i,j=0}^{\infty }$ are the lower and, 
respectively, the upper parts of $S$, such that $L_{ij}=S_{ij}$
for $i>j$ and $L_{ij}=0$ for $i\leq j$, while
$Q_{ij}=S_{ij}$
for $i\leq j$ and $Q_{ij}=0$ for $i>j$. Since, formally
(that is, computing the entries of the arrays that are involved
according to the usual rule of multiplication of operators),
$LC_k=C_kL$ and $QC_k^*=C_k^*Q$
for all $k=1,\ldots,N$, we deduce that
$$\begin{array}{lll}
S-\sum _{k=1}^NC_kSC_k^*&=&L+Q-\sum _{k=1}^NC_k(L+Q)C_k^* \\
 & & \\
 &=&L+Q-\sum _{k=1}^NC_kLC_k^*-\sum _{k=1}^NC_kQC_k^* \\
 & & \\
 &=&L+Q-\sum _{k=1}^NLC_kC_k^*-\sum _{k=1}^NC_kC_k^*Q.
\end{array}$$
But $I=\sum _{k=1}^NC_kC_k^*=P_{\cH_0}$, the orthogonal projection
on $\cH_0$, so we can deduce that
$$\begin{array}{lll}
S-\sum _{k=1}^NC_kSC_k^*&=&L+Q-L(I-P_{\cH_0})-(I-P_{\cH_0})Q \\
 & & \\
 &=&LP_{\cH_0}+P_{\cH_0}Q=GJG^*. 
\end{array}$$

\end{proof}

We now consider truncations of a 
positive definite multi-Toeplitz
kernel on $\FF _N^+$. If $K$ is such a kernel
and $S_{ij}$ are associated to $K$ by the formula \eqref{fourier},
then we define $K_n=[S_{ij}]_{i,j=0}^n$. These are positive
operators on $\oplus _{l=0}^n\cH_l$ 
and their entries satisfy \eqref{condition1}. We also introduce the 
operators $C_{k,n}=P_{\oplus _{l=0}^n\cH_l}C_k\mid \oplus _{l=0}^n\cH_l$.
We obtain another example of a displacement equation.
 
\begin{corollary}\label{C:2.3}
The positive operator $K_n$ satisfies the displacement equation
\begin{equation}\label{basic2}
K_n-\sum _{k=1}^NC_{k,n}K_nC_{k,n}^*=GJG^*,
\end{equation}
where the generator $G$ is given by the formula:
$$G=\left[\begin{array}{cc}
1 & 0 \\
S^*_{01} & S^*_{01} \\
\vdots & \vdots \\
S^*_{0n} &  S^*_{0n} 
\end{array}\right].$$
\end{corollary}

\begin{proof}
This is an immediate consequence of Theorem~\ref{T:2.2}.
\end{proof}

It turns out that the operators satisfying \eqref{basic2} have a 
certain explicit description. For $x\in \oplus _{l=0}^n\cH_l$
we denote by $U(x)$ the upper triangular operator on $\oplus _{l=0}^n\cH_l$
satisfying \eqref{condition1} and whose first row is $x$. 

\begin{theorem}\label{T:2.4}
A positive operator $A$ on $\oplus _{l=0}^n\cH_l$ satisfies \eqref{basic2}
if and only if $A=U(x)^*U(x)-U(y)^*U(y)$
for two vectors $x,y\in \oplus _{l=0}^n\cH_l$
with the property that there is another vector $z\in \oplus _{l=0}^n\cH_l$
such that $U(y)=U(z)U(x)$ and $U(z)$ is a contraction.
\end{theorem}

\begin{proof}
If $A=U(x)^*U(x)-U(y)^*U(y)$ then we use the fact that
for any $x\in \oplus _{l=0}^n\cH_l$ and $k=1,\ldots, N$,
$C_{n,k}U(x)^*=U(x)^*C_{n,k}$, and deduce \eqref{basic2} as
in the proof of Theorem~\ref{T:2.2}.
Conversely, if $A$ satisfies \eqref{basic2}, then we use the fact
that $C_{\sigma ,n}=0$ for $|\sigma |>n$ in order to deduce
that
$$A=GJG^*+\sum _{|\sigma |=1}^nC_{\sigma ,n}GJG^*C_{\sigma ,n}^*.$$
Hence, if $G=[x^*,y^*]$ for some $x,y\in \oplus _{l=0}^n\cH_l$,
then $A=U(x)^*U(x)-U(y)^*U(y)$. Since $A$ is positive, it follows that
there exists a contraction $T$ such that $U(y)=TU(x)$. We notice that
we can assume, without loss of generality, that $x_0\ne 0$ and then
we deduce that $C_{n,k}T^*=T^*C_{n,k}$ for $k=1,\ldots , N$, 
which implies that $T=U(z)$ for some $z\in \oplus _{l=0}^n\cH_l$.
\end{proof}

We mention a simple but useful consequence of this result. 
Let $K$ be a 
positive definite multi-Toeplitz
kernel on $\FF _N^+$ with $K(\emptyset ,\emptyset )=1$
and let $S_{ij}$, $i,j\geq 0$, be the operators associated to $K$
by \eqref{fourier}. Let 
$$x_n=\left[\begin{array}{c}
1  \\
S^*_{01}  \\
\vdots \\
S^*_{0n}  
\end{array}\right], \quad y_n=\left[\begin{array}{c}
 0 \\
 S^*_{01} \\
 \vdots \\
  S^*_{0n} 
\end{array}\right],
$$
then $K_n=U(x_n)^*U(x_n)-U(y_n)^*U(y_n)$ and by the previous result
there exists $z_n\in \oplus _{l=0}^n\cH_l$ such that 
$
U(y_n)=U(z_n)U(x_n).
$
In fact, $z_n$ is uniquely determined and also 
$P_{\oplus _{l=0}^n\cH_l}U(z_{n+1})\mid \oplus _{l=0}^n\cH_l=U(z_n)$ 
for all $n\geq 1$. It follows that there exists a unique
$T(K)\in \cS (\cH)$ such that $P_{\oplus _{l=0}^n\cH_l}T\mid 
\oplus _{l=0}^n\cH_l=U(z_n)$ and $T(K)_{00}=0$.
Denote by $\cS _0(\cH)$ 
the set of those elements $T\in \cS (\cH)$
with $T_{00}=0$.

\begin{corollary}\label{C:2.5}
The mapping $K\rightarrow T(K)$ gives a one-to-one correspondence 
between the set 
of positive definite multi-Toeplitz
kernels on $\FF _N^+$ with $K(\emptyset ,\emptyset )=1$ and $\cS _0(\cH)$.
\end{corollary}

\subsubsection{Pick kernels.}
This type of kernels were recently introduced 
in \cite{Ar}, \cite{Po2}, and they were shown to
have a certain universality property with respect
to a Nevanlinna-Pick problem, \cite{AM}. 

Let ${\cB }^N$ be the open unit ball in ${\CC }^N$.
For distinct points $\lambda _1$, $\ldots $,
$\lambda _L$ in ${\cB }^N$ and complex numbers
$b_1$, $\ldots $, $b_L$, we define the Pick kernel
by the formula
$$R=\left[\frac{1-b_j\overline{b_l}}{1-\langle \lambda _j,\lambda _l
\rangle }\right]_{j,l=1}^L,$$
where $\langle \lambda _j,\lambda _l\rangle $
denotes the Euclidean inner product in ${\CC }^N$.

We can define 
$$F_j=\mbox{diag}(\lambda _l^{(j)})_{l=1}^L,\quad j=1,\ldots ,N,$$
the diagonal matrix with the diagonal made of the $j$th components 
of the points $\lambda _1$, $\ldots $,
$\lambda _L$. Also, we set
$$G=\left[\begin{array}{cc}
1 & b_1 \\
1 & b_2 \\
\vdots & \vdots \\
1 & b_L 
\end{array}
\right].$$
We can check by direct computation that 
\begin{equation}\label{pick}
R-\sum _{k=1}^NF_kRF_k^*=GJG^*.
\end{equation}
It was shown in \cite{CSK} how to use this equation in order to 
solve some multidimensional Nevanlinna-Pick type problems as those in 
\cite{Ar}, \cite{Po2}.

\section{Scattering Experiments}

One of the main results in displacement structure theory
refers to the possibility of associating a scattering  
operator to the data given by the generators of the equation 
\eqref{displacement}. In this section we prove a similar result
for displacement equations of the type of \eqref{pick} 
and we explain the connection with a special class of time variant
linear systems.

Let $\cE _1$, $\cE _2$, and $\cG $ be Hilbert spaces.  
We consider a displacement equation of the form 
\begin{equation}\label{newdisplacement}
A-\sum _{k=1}^NF_kAF_k^*=GJG^*,
\end{equation}
where $F_k\in \cL(\cG)$, $k=1,\ldots ,N$, are given contractions
on the Hilbert
space $\cG$. Also $G=\left[\begin{array}{cc}
U & V
\end{array}\right]\in \cL (\cE _1\oplus \cE _2,\cG)$
and $J=I_{\cE _1}\oplus I_{\cE _2}$.
The wave operators associated to \eqref{newdisplacement} 
are introduced by the formulae:
$U^*_{\infty }=\left[U_k\right]_{k=0}^{\infty }$,
$V^*_{\infty }=\left[V_k\right]_{k=0}^{\infty }$, where
$U_k=\left[F_{\sigma }U\right]_{|\sigma |=k}:
\cH _k\otimes \cE _1\rightarrow \cG$,  
$V_k=\left[F_{\sigma }V\right]_{|\sigma |=k}:
\cH _k\otimes \cE _2\rightarrow \cG$.
We will assume that both $U_{\infty } $ and $V_{\infty } $ are bounded
and also that $\lim _{k\rightarrow \infty }
\sum _{|\sigma |=k}\|F^*_{\sigma }g\|=0$
for all $g\in \cG $. Under these assumptions we
deduce that
\eqref{newdisplacement} has a unique solution given by
\begin{equation}\label{solution}
A=U^*_{\infty }U_{\infty }-V^*_{\infty }V_{\infty }. 
\end{equation}

The next result was already noticed in \cite{CSK}, but here
we present a different approach based on system theoretic ideas.

\begin{theorem}\label{T:3.1}
The solution \eqref{solution}
of the displacement equation \eqref{newdisplacement} 
is positive if and only if there exists 
$T\in \cS (\cH ,\cE _1,\cE _2)$ such that $V_{\infty }=TU_{\infty }$.
\end{theorem}

\begin{proof}
Assume $A=U^*_{\infty }U_{\infty }-V^*_{\infty }V_{\infty }\geq 0$ 
and let $A=LL^*$ be 
a factorization of $A$ with $L\in \cL (\cF ,\cG )$ for some 
Hilbert space $\cF $. From \eqref{newdisplacement} we deduce
that 
$$LL^*+VV^*=\sum _{k=1}^NF_kLL^*F^*_k+UU^*.$$
In matrix form, 
\begin{equation}\label{base}
\left[\begin{array}{cc}
L & V
\end{array}\right]
\left[\begin{array}{c}
L^* \\
V^*
\end{array}\right]=
\left[\begin{array}{cccc}
F_1L & \ldots & F_NL & U
\end{array}\right]
\left[\begin{array}{c}
L^*F^*_1 \\
\vdots \\
L^*F^*_N \\
 U^*
\end{array}\right].
\end{equation}  
Define $A^*=\left[\begin{array}{cc}
L & V
\end{array}\right]$ and 
$B^*=\left[\begin{array}{cccc}
F_1L & \ldots & F_NL & U
\end{array}\right]$, then we deduce from \eqref{base}
that there exists an unitary operator
$\theta _0\in \cL (\overline{\cR (B)},\overline{\cR (A)})$ such that
$A=\theta _0B.$
It follows that there exist Hilbert spaces 
$\cR _1$, $\cR _2$, and an unitary extension 
$\theta \in \cL (\cF ^{\oplus N}\oplus \cE _1 \oplus \cR _1,
\cF \oplus \cE _2 \oplus \cR _2)$ of $\theta _0$, hence this extension  
satisfies the relation 
\begin{equation}\label{relation}
\left[\begin{array}{c}
A \\
0_{\cR _2}
\end{array}\right]
=\theta \left[\begin{array}{c}
B \\
0_{\cR _1}
\end{array}\right].
\end{equation}

\noindent
Let $\theta _{ij}$, $i\in \{1,2,3\}$, $j\in \{1,2,\ldots ,N+2\}$,
be the matrix coefficients of $\theta $. It is convenient
to rename some of these coefficients. Thus, we set
$$\begin{array}{cccc}
X_k=\theta _{1k}, & k=1,\ldots N, & \quad & Z=\theta _{1,N+1}, \\
 & &  & \\
Y_k=\theta _{2k}, & k=1,\ldots N, & \quad & W=\theta _{2,N+1}.
\end{array}$$
From \eqref{relation} we deduce that
$$L^*=\sum _{k=1}^NX_kL^*F^*_k+ZU^*$$
and 
$$V^*=\sum _{k=1}^NY_kL^*F^*_k+WV^*.$$
By induction we deduce that 
\begin{equation}\label{key}
V^*=WU^*+\sum _{k=1}^N\sum _{|\sigma |=0}^nY_kX_{\sigma }ZU^*F_{k\sigma }^*
+\sum _{|\tau |=n+1}Q_{\tau }L^*F_{\tau }^*,
\end{equation}
where $Q_{\tau }$ are monomials of length $|\tau |$ in the variables
$X_1$, $\ldots $,$X_N$, $Y_1$, $\ldots $, $Y_N$.
Since $\theta $ is unitary it follows that all $Q_{\tau }$ are contractions.

We define $T_{00}=W$ and for $j>0$,
$$T_{0j}=[Y_kX_{\sigma }Z]_{|\sigma |=j-1 ;k=1,\ldots ,N}.$$
Then we define $T_{ij}$, $i>0$, $j\geq i$, by the 
formula \eqref{condition} and $T_{ij}=0$ for $i>j$. We show that 
$T=[T_{ij}]_{i,j=0}^{\infty }$ belongs to $\cS (\cH ,\cE _1,\cE _2)$. 
To that end, we introduce the notation: for $n \geq 0$,
$$\begin{array}{lcl}
X_1(n)=\left[\begin{array}{ccc}
X_1 & \ldots & X_N
\end{array}\right]^{\oplus N^n},
 & \quad &  Y_1(n)=\left[\begin{array}{ccc}
Y_1 & \ldots & Y_N
\end{array}\right]^{\oplus N^n}, \\
 & & \\
X_2(n)=Z^{\oplus N^n}, & \quad & Y_2(n)=W^{\oplus N^n}.
\end{array}$$
We check by induction that $T_{00}=Y_2(0),$
$T_{01}=Y_1(0)X_2(1)$, and for $j\geq 2$, 
$$T_{0j}=Y_1(0)X_1(1)\ldots X_1(j-1)X_2(j).$$
Consequently, $T^*$ is the transfer map of the linear time variant 
system
$$\left\{\begin{array}{rcl}
x(n+1)&=&X_1(n)^*x(n)+Y_1(n)^*u(n) \\
 & & \\
y(n)&=&X_2(n)^*x(n)+Y_2(n)^*u(n), \quad n\geq 0.
\end{array}\right.
$$
Since each matrix $\left[\begin{array}{cc}
X_1(n) & X_2(n) \\
Y_1(n) & Y_2(n) 
\end{array}
\right]$, $n\geq 0$, is a contraction, it follows that 
$T^*$ is a contraction, hence $T$ is a contraction and 
$T\in \cS (\cH ,\cE _1,\cE _2)$.
Also, since 
$\lim _{k\rightarrow \infty }
\sum _{|\sigma |=k}\|F^*_{\sigma }g\|=0$
for all $g\in \cG $, we deduce from \eqref{key}
that $V_{\infty }=TU_{\infty }$.
\end{proof}

This result is a generalization of Theorem~\ref{T:3.1}. It 
can be used to solve interpolation and moment problems in several 
non commuting
variables as suggested in \cite{CSK}.
The proof given here emphasizes the fact that this type of
problems appears as an interesting particular case
of similar
questions in time variant setting (as described, for instance, in 
\cite{Co}). Next we develop this idea by giving a
parametrization of the elements in $\cS (\cH ,\cE _1,\cE _2)$.

\subsection{Schur Parameters}
We introduce a Schur type parametrization for   
$\cS(\cH,\cE _1,\cE _2)$ and then give a number of applications. 
A first instance of such a result might be considered
Euler's description of $SO(3)$. Schur's classical result
in \cite{Sc} gives a parametrization of the contractive
holomorphic functions on the unit disk and Szeg\"o's theory 
of orthogonal polynomials provides an alternative to this result
for probability measures on the unit circle.

Here we provide a refinement of a parametrization
of $\cS (\cH ,\cE _1,\cE _2)$ given in \cite{Po1}.
We can deal with a slightly more general situation than 
$\cS (\cH )$. For a contraction 
$C$ we denote by $\cD _C$ its defect space, that is
the closure of the range of the defect operator $D_C=(I-C^*C)^{1/2}$.
Let $\Pi (\cE _1,\cE _2)$ be the family of sequences 
$\{\gamma _{\sigma }\}_{\sigma \in \FF _N^+}$ with the 
properties: $\gamma _{\emptyset }\in \cL (\cE _1,\cE _2)$ and for $|\sigma |>0$,  
\begin{equation}\label{compatibility}
\gamma _{k\sigma }\in \cL (\cD _{\gamma _{\sigma }}, 
\cD _{\gamma ^*_{k\sigma -1}}),
\end{equation}
where $k=1,\ldots N$ and 
$k\sigma -1$ denotes the predecessor of $k\sigma $
with respect to the lexicographic order.

\begin{theorem}\label{T:3.2}
There exists a one-to-one correspondence between 
$\cS (\cH ,\cE _1,\cE _2)$ and 
$\Pi (\cE _1,\cE _2)$.
\end{theorem}

\begin{proof}
A simple proof can be obtained as an application
of the Schur parametrization of upper triangular contractions
in \cite{Co}.
Let $T\in \cS(\cH ,\cE _1,\cE _2)$ and for each $j\geq 0$, $T_{0j}$ 
is the row matrix $[c_{\sigma }]_{|\sigma |=j}$.
We use Theorem~2.2.1 in \cite{Co} in order to associate
to $T$ a family $\Gamma $ of Schur parameters. The Schur parameters
of the contraction
$$\left[\begin{array}{cc}
T_{00} & T_{01} \\
0 & T_{11} 
\end{array}\right]=\left[\begin{array}{ccccc}
c_{\emptyset } & c_1 & \ldots & & c_N \\
0 & c_{\emptyset }& 0 & \ldots & 0 \\
 & \ddots & \ddots & \ddots & \\
0 & \ldots & & 0 & c_{\emptyset }
\end{array}\right]$$
are given by $\Gamma _{00}=\gamma _{\emptyset }=c_{\emptyset }$,
$\Gamma _{11}=\Gamma _{00}^{\oplus N}$ and 
$\Gamma _{01}\in \cL (\cD ^{\oplus N}_{\gamma _{\emptyset }},
\cD _{\gamma _{\emptyset }})$. 

For the next step
we associate Schur parameters to the contraction
$$\left[\begin{array}{ccc}
T_{00} & T_{01} & T_{02} \\
0 & T_{11} & T_{12} \\
0 & 0 & T_{22}
\end{array}\right].
$$
We notice that 
$$
\left[\begin{array}{cc}
T_{11} & T_{12} \\
0 & T_{22} 
\end{array}\right]
=\left[\begin{array}{cc}
T_{00}^{\oplus N} & T_{01}^{\oplus N} \\
0 & T_{11}^{\oplus N} 
\end{array}\right].
$$
By the formula (2.2.3) in \cite{Co}, the entries 
of the matrix 
$\left[\begin{array}{cc}
T_{11} & T_{12} \\
0 & T_{22} 
\end{array}\right]$ are polynomials in the Schur parameters, their
adjoints and their defect operators. Consequently, 
the contraction associated to the Schur parameters
$\Gamma _{00}^{\oplus N}$, $\Gamma _{01}^{\oplus N}$ and 
$\Gamma _{11}^{\oplus N}$ is precisely
$\left[\begin{array}{cc}
T_{00}^{\oplus N} & T_{01}^{\oplus N} \\
0 & T_{11}^{\oplus N} 
\end{array}\right].$
From the uniqueness of the Schur parameters, we deduce 
that the Schur parameters of 
$\left[\begin{array}{cc}
T_{11} & T_{12} \\
0 & T_{22} 
\end{array}\right]$
are precisely 
$\Gamma _{00}^{\oplus N}$, $\Gamma _{01}^{\oplus N}$ and 
$\Gamma _{11}^{\oplus N}$. By induction, we deduce that we can associate
to $T$ the family of Schur parameters $\Gamma =\{\Gamma _{ij}\}_{i\leq j}$
such that each $\Gamma _{ij}$ satisfies \eqref{condition}.
We notice that $\Gamma _{0j}:\cD_{\Gamma _{0,j-1}}^{\oplus N}\rightarrow
\cD _{\Gamma ^*_{0,j-1}}$ for $j\geq 1$.
Also, each $\Gamma _{0j}$ is a row contraction,
$\Gamma _{0j}=[\gamma ^0_{\sigma }]_{|\sigma |=j}$.
By Proposition~1.4.2 in \cite{Co}, 
$$\Gamma _{01}=\left[\begin{array}{cccc}
\gamma _1 & D_{\gamma _1^*}\gamma _2 & \ldots  &
D_{\gamma _1^*}\ldots D_{\gamma _{N-1}^*}\gamma _N
\end{array}\right],
$$
where $\gamma _k\in \cL (\cD _{\gamma _{\emptyset }},\cD _{\gamma ^*_{k-1}})$,
$k=1,\ldots ,N$. Also, there exist unitary operators 
$u_1:\cD _{\Gamma _{01}}\rightarrow \oplus _{k=1}^N\cD _{\gamma _k}$ and
$v_1:\cD _{\Gamma ^*_{01}}\rightarrow \cD _{\gamma ^*_N}$.
Next, we have that $\Gamma _{02}:\cD _{\Gamma _{01}}^{\oplus N}
\rightarrow \cD _{\Gamma ^*_{01}}$ and we can define the row
contraction $\tilde \Gamma _{02}=[\gamma ^0_{\sigma }]_{|\sigma |=2}$
from $\left( \oplus _{k=1}^N\cD _{\gamma _k}\right) ^{\oplus N}$
into $\cD _{\gamma ^*_N}$ by the formula $\tilde \Gamma _{02}=
v_1\Gamma _{02}(u_1^{\oplus N})^*$.
As above, 
$$\tilde \Gamma _{02}=\left[
\begin{array}{cccc}
\gamma _{11} & D_{\gamma ^*_{11}}\gamma _{12} & \ldots &
D_{\gamma ^*_{11}}\ldots D_{\gamma ^*_{N,N-1}}\gamma _{NN}
\end{array}\right],
$$
where $\gamma _{k\tau }
\in \cL (\cD_{\gamma _{\tau }},\cD _{\gamma ^*_{k\tau -1}})$,
$k=1, \ldots ,N$, $|\tau |=1$.
We proceed by induction in order to construct the whole family
of Schur parameters $\{\gamma _{\sigma }\}_{\sigma \in \FF _N^+}$ 
with $\gamma _{\emptyset }\in \cL (\cE _1,\cE _2)$ and obeying
\eqref{compatibility}.
\end{proof}

The parameters
$\{\Gamma _{0j}\}_{j\geq 0}$ are the Schur parameters
associated to $T$ in \cite{Po1}. The above proof explain 
their construction within the framework of \cite{Co}.
The new parameters $\gamma _{\sigma }$ are more like "true"
Schur parameters in the sense that they are contractions 
subject to the only restrictions given by 
the compatibility relations \eqref{compatibility}.

A similar parametrization can be obtained for the set of positive 
multi-Toeplitz kernels. Denote by $\Pi _0(\cE _1,\cE _2)$ the set of
families $\{\gamma _{\sigma }\}_{\sigma \in \FF _N^+}$
in $\Pi (\cE _1,\cE _2)$ with $\gamma _{\emptyset }=0$.

\begin{theorem}\label{T:3.3}
There exists a one-to-one correspondence between the set
of $\cE$-valued positive multi-Toeplitz kernels 
with $K(\emptyset, \emptyset )=1$ and 
$\Pi _0(\cE )$.
\end{theorem} 

\begin{proof}
This is a consequence of Theorem~\ref{T:3.2} and Corollary~\ref{C:2.5}.
A direct proof can be obtained as an application of 
Theorem~1.5.3 in \cite{Co}. The parametrization given here
is a refinement of one given in \cite{Po2}. 
\end{proof}

An application gives another representation of the 
class $\cS (\cH ,\cE _1,\cE _2)$.
Motivated by the proof of Theorem~\ref{T:3.1}
we consider the class $\cB \cS _N(\cE _1,\cE _2,\cF )$ of systems
$$\Sigma \left\{\begin{array}{rcl}
x(n+1)&=&A(n)x(n)+B(n)u(n) \\
 & & \\
y(n)&=&C(n)x(n)+D(n)u(n), \quad n\geq 0,
\end{array}\right.
$$
where $\cE _1$, $\cE _2$, $\cF $ are Hilbert spaces,
$A(0)\in \cL (\cF ,\cF ^{\oplus N})$,
$B(0)\in \cL (\cE_2 ,\cF ^{\oplus N})$, 
$C(0)\in \cL (\cF ,\cE_1)$, $D\in \cL (\cE_2, \cE _1)$,
and $A(n)=A(n-1)^{\oplus N}$,
$B(n)=B(n-1)^{\oplus N}$,
$C(n)=C(n-1)^{\oplus N}$,
$D(n)=D(n-1)^{\oplus N}$ (these relations are all of the type of
\eqref{condition}).
Also we assume that for all $n\geq 0$, the operator
$\left[\begin{array}{cc}
A(n) & B(n) \\
C(n) & D(n) 
\end{array}\right]$
is a coisometry. The evolution of the system
$\Sigma $ can be encoded by the transfer map $T_{\Sigma }$ that shows
how the infinite input sequence $(u(0),u(1),\ldots )$ is transformed
by the system into an output sequence $(y(0),y(1),\ldots )$ provided
that $x(0)=0$. It is easy to see that the transfer map 
$T_{\Sigma }$ is a lower triangular contraction with the matrix
coefficients $T_{\Sigma }(n,n)=D(n)$, $T_{\Sigma }(n+1,n)=
C(n+1)B(n)$ and for $j>1$, 
$T_{\Sigma }(n+j,n)=C(n+j)A(n+j-1)\ldots A(n+1)B(n)$.

\begin{theorem}\label{T:3.4}
$T$ belongs to $\cS (\cH ,\cE _1,\cE _2)$ if and only if 
$T^*$ is the transfer map of a system in $\cB \cS _N(\cE _1,\cE _2,\cF)$
for some Hilbert space $\cF $.
\end{theorem}

\begin{proof}
We use the construction from the proof
of Theorem~2.3.3 in \cite{Co}. Thus, 
let $\{\Gamma _{ij}\}_{i\leq j}$ be the Schur parameters
associated to $T$ in the proof of Theorem~\ref{T:3.2}.
For each $i\geq 0$ we consider the row contraction
$$L_i=\left[\begin{array}{cccc}
\Gamma _{ii} & D_{\Gamma ^*_{ii}}\Gamma _{i,i+1} &
D_{\Gamma ^*_{ii}}D_{\Gamma ^*_{i,i+1}}\Gamma _{i,i+2} & \ldots
\end{array}\right].
$$
Let $V_i$ be the isometry associated to $L_i$ by the 
formula 1.6.3 in \cite{Co}, that is,
\begin{equation}\label{v}
V_i=
\left[\begin{array}{cccc}
\Gamma _{ii} & D_{\Gamma ^*_{ii}}\Gamma _{i,i+1} &
D_{\Gamma ^*_{ii}}D_{\Gamma ^*_{i,i+1}}\Gamma _{i,i+2} & \ldots \\
D_{\Gamma _{ii}} & -\Gamma ^*_{ii}\Gamma _{i,i+1} & 
-\Gamma ^*_{ii}D_{\Gamma ^*_{i,i+1}}\Gamma _{i,i+2} & \ldots \\
0 & D_{\Gamma _{i,i+1}} & -\Gamma ^*_{i,i+1}\Gamma _{i,i+2} & \ldots \\
0 & 0 & D_{\Gamma _{i,i+2}} & \\
\vdots & \vdots & & \ddots 
\end{array}\right]=
\left[\begin{array}{cc}
\Gamma _{ii} & B(i)^* \\
C(i)^* & A(i)^*
\end{array}\right].
\end{equation}
$V_0$ is defined on the space $\cE \oplus \oplus _{k=0}^{\infty }
\cD_{\Gamma _{0k}}^{\oplus N}$ with values in 
in $\cE \oplus \oplus _{k=0}^{\infty }
\cD_{\Gamma _{0k}}$. We set $\cF =\oplus _{k=0}^{\infty }
\cD_{\Gamma _{0k}}$ and let $\phi :(\oplus _{k=0}^{\infty }
\cD_{\Gamma _{0k}})^{\oplus N}\rightarrow 
\oplus _{k=0}^{\infty }
\cD_{\Gamma _{0k}}^{\oplus N}$ be the unitary operator defined 
by the formula
$\phi (\oplus _{l=1}^N\oplus _{k=0}^{\infty }u_{kl})=
\oplus _{k=0}^{\infty }\oplus _{l=1}^Nu_{kl}$.
Then we define 
$$
\left[\begin{array}{cc}
\bar D^*(i)& \bar B^*(i) \\
\bar C^*(i) & \bar A^*(i)
\end{array}\right]=
\left[\begin{array}{cc}
\Gamma _{ii}& B^*(i)\phi ^{\oplus N^i} \\
(\phi ^{\oplus N^{i-1}})^*C^*(i) & 
(\phi ^{\oplus N^{i-1}})^*A^*(i)\phi ^{\oplus N^i}
\end{array}\right].
$$
By Theorem~2.3.3 in \cite{Co},
$T^*_{i,i+1}=C(i+1)B(i)=\bar C(i+1)\bar B(i)$ and for $j>i$,
\begin{equation}\label{transfer}
T^*_{ij}=C(j)A(j-1)\ldots A(i+1)B(i)=
\bar C(j)\bar A(j-1)\ldots \bar A(i+1)\bar B(i).
\end{equation}
Since the elements of the family $\Gamma $
satisfy \eqref{condition} and the matrix coefficients
of $V_i$ are monomials in $\Gamma _{ij}$, their adjoints and
defect operators, it follows that the system associated
to the families $\{\bar A(n)\}_{n\geq 0}$,
$\{\bar B(n)\}_{n\geq 0}$,
$\{\bar C(n)\}_{n\geq 0}$ and 
$\{\bar D(n)\}_{n\geq 0}$ belongs to $\cB \cS (\cE ,\cF)$
and $T^*$ is the transfer map of this system.
\end{proof}

\section{Recursive Constructions}

As noticed in the proof of Theorem~\ref{T:3.4}, the Schur 
parameters introduced in Theorem~\ref{T:3.2} are basically given by the 
parametrization of the unitary group. A convenient way to encode
the dependence of the matrix coefficients of $T$ on the Schur paramaters
is given by a continued fraction algorithm introduce for $N=1$
by Schur in \cite{Sc}. Here we describe such an algorithm for $N>1$.
Of course, the main dificulty consists in keeping all the 
contructions within the tensor algebra. Due to the fact
that there is a standard filtration on the tensor algebra, 
the contructions do not go through as in the one dimensional case, still
we can obtain a sort of graded Schur type algorithm.
For the algebra of upper triangular matrices, a Schur type
algorithm was known (as a consequence of a similar
algorithm for the spectral factorization of a family of matrices with
displacement structure, see \cite{KS} or \cite{Co} for a detailed 
discussion). If we use it directly 
to our situation, it will produce the paramaters $\{\gamma _{\sigma }\}_
{\sigma \in \FF _N^+}$ of Theorem~\ref{T:3.2}, but the objects 
produced at each step are no longer in the tensor algebra. In this section
we show how to remedy this situation.

We first notice another useful representation of a 
transfer map (we will assume in this section that $\cE _1=\cE _1=\CC $). 
It is convenient to extend the symbol map \eqref{izom}. 
Thus, for $T\in \cL (\oplus _{k=0}^{\infty }\cH _k\otimes \cE _1,
\oplus _{k=0}^{\infty }\cH _k\otimes \cE _2 )$ an upper triangular 
operator satisfying \eqref{condition}, we write
$$\Phi (\left[T_{0j}\right]_{j=0}^{\infty })=T.$$

Now, if $T_{\Sigma }$ is the transfer map 
of a system $\Sigma $ in $\cB \cS _N(\cF )=\cB \cS _N(\CC ,\CC ,\cF )$, then 
$$\begin{array}{rcl}
T^*_{\Sigma }&=&
\Phi (D(0)^*,0,\ldots ) 
+
\Phi (0,B(0)^*,0,\ldots ) \\
& & \\ 
& & \times 
(I-\Phi (0,A(0)^*,0,\ldots ))^{-1}
\Phi (C(0)^*,0,\ldots ),
\end{array}
$$
provided that $\|\Phi (0,A(0)^*,0,\ldots )\|<1$. 
We notice that this is 
the case if and only if $\|A(0)\|<1$. Let 
$A(0)^*=\left[A_j(0)^*\right]_{j=1}^N$, 
$B(0)^*=\left[B_j(0)^*\right]_{j=1}^N$. Then
$$\Phi (0,A(0)^*,0,\ldots )=
\sum _{j=1}^N\Phi (A_j(0)^*,0,\ldots )C_j^-$$
and 
$$\Phi (0,B(0)^*,0,\ldots )=
\sum _{j=1}^N\Phi (B_j(0)^*,0,\ldots )C_j^-,$$
so that 
\begin{equation}\label{realizare}
\begin{array}{rcl}
T^*_{\Sigma }&=&\Phi (D(0)^*,0,\ldots )+
(\sum _{j=1}^N\Phi (B_j(0)^*,0,\ldots )C_j^-) \\
& & \\
& &\times  
(I-\sum _{j=1}^N\Phi (A_j(0)^*,0,\ldots )C_j^-)^{-1}
\Phi (C(0)^*,0,\ldots ),
\end{array}
\end{equation}
provided that $\|A(0)\|<1$.
Next, we prove a refinement of \eqref{realizare}
which is the key step of a Schur type algorithm.
We introduce the translation operator given by the matrix
$$S=\left[\begin{array}{cccc}
0 & 0 &  & \ldots \\ 
1 & 0 &  & \ldots \\
0 & 1 & 0 & \ldots  \\
\vdots & &  \ddots  
\end{array}
\right].
$$

\begin{theorem}\label{T:4.1}
If $T$ belongs to $\cS (\cH )$ then there 
exists an unique $T_1\in \cS (\cH , \cD ^{\oplus N}_{T_{00}},
\cD _{T^*_{00}})$ such that
$$
\begin{array}{rcl}
T&=&\Phi (T_{00},0,\ldots )+
\Phi (D_{T^*_{00}},0,\ldots )T_1S \\
 & & \\
& & \times 
(I-\Phi (-T^*_{00},0,\ldots )T_1S
\Phi (D_{T_{00}},0,\ldots ).
\end{array}
$$
\end{theorem}

\begin{proof}
Let $\{\Gamma _{ij}\}_{i\leq j}$ be the Schur parameters associated 
to $T$ as in the proof of Theorem~\ref{T:3.2}. Then $|\Gamma _{00}|\leq 1$.
If $|\Gamma _{00}|=1$, then we must have $\Gamma _{ij}=0$ for $j>i$, so that
$T=\Phi (\Gamma _{00},0,\ldots )$ and $T_1=\ldots =T_N=0$. 
Assume $|\Gamma _{00}|<1$. Then $T^*$ is the transfer map of the system
associated to the operators $V_i$ defined by 
\eqref{v}. Thus, we have

$$\begin{array}{rcl}
V_0&=&\left[\begin{array}{ccc}
\Gamma _{00} & B(0)^* \\
C(0)^* & A(0)^*
\end{array}\right]=
\left[\begin{array}{cccc}
\Gamma _{00} & D_{\Gamma ^*_{00}}\Gamma _{01} &
D_{\Gamma ^*_{00}}D_{\Gamma ^*_{01}}\Gamma _{02} & \ldots \\
D_{\Gamma _{00}} & -\Gamma ^*_{00}\Gamma _{01} & 
-\Gamma ^*_{00}D_{\Gamma ^*_{01}}\Gamma _{02} & \ldots \\
0 & D_{\Gamma _{01}} & -\Gamma ^*_{01}\Gamma _{02} & \ldots \\
0 & 0 & D_{\Gamma _{02}} &  \\
\vdots & \vdots & & \ddots 
\end{array}\right] \\
 & & \\
&=&
\left[\begin{array}{cccc}
\Gamma _{00} & D_{\Gamma ^*_{00}} & 0 & \ldots  \\
D_{\Gamma _{00}} & -\Gamma ^*_{00} & 0 & \ldots  \\
0 & 0 & I & \\
\vdots & \vdots & & \ddots 
\end{array}\right]
\left[\begin{array}{cccc}
I & 0 & 0 & \ldots  \\
0 & \Gamma _{01} & D_{\Gamma ^*_{01}}\Gamma _{02} & \ldots  \\
0 & D_{\Gamma _{01}} & -\Gamma ^*_{01}\Gamma _{02} &  \\
\vdots & \vdots & & \ddots 
\end{array}\right] \\
 & & \\
&=&
\left[\begin{array}{ccc}
\Gamma _{00} & D_{\Gamma ^*_{00}} & 0 \\
D_{\Gamma _{00}} & -\Gamma ^*_{00} & 0  \\
0 & 0 & I 
\end{array}\right]
\left[\begin{array}{ccc}
I & 0 & 0  \\
0 & \Gamma _{01} & B'(0)^*  \\
0 & C'(0)^* & A'(0)^*
\end{array}\right] \\
&=&\left[\begin{array}{ccc}
\Gamma _{00} & D_{\Gamma ^*_{00}}\Gamma _{01} & D_{\Gamma ^*_{00}}B'(0)^* \\
D_{\Gamma _{00}} & -\Gamma ^*_{00}\Gamma _{01} & -\Gamma ^*_{00}B'(0)^*  \\
0 & C'(0)^* & A'(0)^*  
\end{array}\right].
\end{array}
$$
We assume first that $\|A(0)\|<1$. Then formula \eqref{realizare} 
gives 
$$\begin{array}{rcl}
T&=&\Phi (\Gamma _{00},0,\ldots )+\Phi (0,B(0)^*,0,\ldots ) \\
& & \\ 
& & \times 
\left(I-\Phi (0,A(0)^*,0,\ldots )\right)^{-1}\Phi (C(0)^*,0,\ldots ) \\
 & & \\
&=&
\Phi (\Gamma _{00},0,\ldots )+
\left[\begin{array}{cc}
\Phi (0,D_{\Gamma ^*_{00}}\Gamma _{01},0,\ldots ) & 
\Phi (0,D_{\Gamma ^*_{00}}B'(0)^*,0,\ldots )
\end{array}\right] \\
& & \\ 
& & \times 
\left[\begin{array}{cc}
I-\Phi (0,-\Gamma ^*_{00}\Gamma _{01},0,\ldots ) & 
-\Phi (0,-\Gamma ^*_{00}B'(0)^*,0,\ldots ) \\
-\Phi (0,C'(0)^*,0,\ldots ) & I-\Phi (0,A'(0)^*,0,\ldots )
\end{array}
\right]^{-1}
\left[\begin{array}{c}
\Phi (D_{\Gamma _{00}},0,\ldots ) \\ 
0
\end{array}\right].
\end{array}
$$
Since $\|A(0)\|<1$, we also have $\|A'(0)\|<1$ so that we 
can use a well-known formula for the inversion of a $2\times 2$
matrix in order to deduce that 
$$\begin{array}{rcl}
T&=&\Phi (\Gamma _{00},0,\ldots )+\Phi (D_{\Gamma ^*_{00}},0,\ldots ) 
(\Phi (0,\Gamma _{00},0,\ldots )\Delta ^{-1} \\
 & & \\ 
& &+ 
\Phi (0,B'(0)^*,0,\ldots )(I-\Phi (0,A'(0)^*,0,\ldots ))^{-1} \\
 & & \\
 & & \times 
\Phi (0,C'(0)^*,0,\ldots )\Delta ^{-1})
\Phi (D_{\Gamma _{00}},0,\ldots )
\end{array}
$$
where 
$$\begin{array}{rcl}
\Delta &=&
I-\Phi (0,-\Gamma ^*_{00}\Gamma _{01},0,\ldots )-
\Phi (0,-\Gamma ^*_{00}B'(0)^*,0,\ldots ) \\
 & & \\
 & & \times 
(I-\Phi (0,A'(0)^*,0,\ldots ))^{-1}\Phi (0,C'(0)^*,0,\ldots ) \\
 & & \\
&=&I-\Phi (0,-\Gamma ^*_{00},0,\ldots )
(\Phi (0,-\Gamma _{00},0,\ldots ) \\
& & \\
& &+ 
\Phi (0,B'(0)^*,0,\ldots ) 
(I-\Phi (0,A'(0)^*,0,\ldots ))^{-1}\Phi (0,C'(0)^*,0,\ldots ))
\end{array}
$$
Since $\|\Phi (0,-\Gamma ^*_{00},0,\ldots )\|<1$, $\Delta $ is invertible
and the previous calculation makes sense. Also, we define
$$R=\Phi (0,\Gamma _{01},0,\ldots )+
\Phi (0,B'(0)^*,0,\ldots ) 
(I-\Phi (0,A'(0)^*,0,\ldots ))^{-1}\Phi (0,C'(0)^*,0,\ldots ).$$
Since 
$\left[\begin{array}{cc}
\Gamma _{01} & B'(0)^* \\
C'(0)^* & A'(0)^*
\end{array}
\right]$ is a contraction (in fact, it is an isometry), we can easily
check that $R\in \cS _0(\cH )$. Actually, we can notice that $R$ is 
obtained from the Schur parameters $\{\tilde{\Gamma }_{ij}\}_{i\geq j}$, 
where $\tilde{\Gamma }_{ii}=0$ and $\tilde{\Gamma }_{ij}=\Gamma _{ij}$
for $j>i$. We obtain that
\begin{equation}\label{schur}
T=\Phi (\Gamma _{00},0,\ldots )
+\Phi (D_{\Gamma ^*_{00}},0,\ldots )
R(I-\Phi (-\Gamma ^*_{00},0,\ldots )R)^{-1}
\Phi (D_{\Gamma _{00}},0,\ldots ).
\end{equation}
Since $\|\Phi (-\Gamma ^*_{00},0,\ldots )\|<1$, the previous 
formula makes sense without our assumption that $\|A(0)\|<1$. In fact, 
an approximation argument shows that \eqref{schur} holds without
that assumption. 
Thus, let $\rho \in (0,1)$ and consider $T^*_{\rho }$ the transfer 
map of the system $\{\rho A(n),B(n),C(n),D(n)\}$ (not necessarely 
in $\cB\cS_N(\cH )$). All the previous calculations go through
and we obtain a contraction $R_{\rho }$ such that \eqref{schur}
holds for $T_{\rho }$. Letting $\rho \rightarrow 1$ we obtain 
\eqref{schur}. 
\end{proof} 

This result leads to the following Schur type algorithm 
for a $T\in \cS (\cH)$.
$$
\left\{\begin{array}{l}
 T_0=T, \quad \gamma _n=T_{n,00}, \,\, n\geq 0, \\
 \\
 T_{n+1}=\mbox{the unique solution of the equation} \\
\quad \quad \quad \quad \quad \quad \quad \quad   T_n=\Phi (\gamma _n,0,\ldots ) \\
\quad +\Phi (D_{\gamma ^*_n},0,\ldots )T_{n+1}S^{\oplus N^n}
(I-\Phi (-\gamma ^*_n,0,\ldots )T_{n+1}S^{\oplus N^n})^{-1}
\Phi (D_{\gamma _n},0,\ldots ).
\end{array}
\right.
$$

The algorithm generates a sequence of elements
$T_n\in \cS (\cH ,\cD _{T^{\oplus N}_{n-1}}, \cD _{T^*_{n-1}})$
and it is easily seen that $\gamma _0=\Gamma _{0n}$, $n\geq 0$.
Therefore $\{\gamma _n\}_{n\geq 0}$
uniquely determine $T$ and the algorithm encodes the dependence
of $T$ on these parameters.





\end{document}